\documentclass[20pt,twoside,draft
]{article}
\usepackage{cmap}                        % ГѓВЏГѓВ®Гѓв‚¬Гѓв‚¬ГѓВҐГѓВ°ГѓЕ ГѓВЄГѓВ  ГѓВЇГѓВ®ГѓЕЎГѓВ±ГѓВЄГѓВ  ГѓВ°ГѓВіГѓВ±ГѓВ±ГѓВЄГѓЕЎГѓВµ ГѓВ±ГѓВ«ГѓВ®ГѓВў ГѓВў PDF (pdflatex)
\usepackage[cp1251]{inputenc}            % ГѓВ‚ГѓВ»ГѓВЎГѓВ®ГѓВ° ГѓВїГѓВ§ГѓВ»ГѓВЄГѓВ  ГѓЕЎ ГѓВЄГѓВ®Гѓв‚¬ГѓЕЎГѓВ°ГѓВ®ГѓВўГѓВЄГѓЕЎ
\usepackage[english, russian]{babel}
\usepackage[left=1cm,right=1cm,top=1cm,bottom=2cm]{geometry} % ГѓВЇГѓВ®ГѓВ«ГѓВї ГѓВ±ГѓВІГѓВ°ГѓВ ГѓВ­ГѓЕЎГѓВ¶ГѓВ»
\usepackage[russian]{babel}
\usepackage{indentfirst, array}
\usepackage{amscd,latexsym}

\usepackage{tabularx}
\usepackage{multirow}
\usepackage{amsmath,amsfonts,amssymb,amsthm,indentfirst,enumerate,textcomp}
\usepackage{graphicx}
\addto\captionsrussian{}

\theoremstyle{plain}
\newtheorem{theorem}{Theorem}

\begin{document}

\begin{center}
\begin{Large}
\textbf{Recent advances in Bergman type projections in bounded pseudoconvex and tubular domains.}\\
\end{Large}
\vskip 0.3cm
\textbf{R.F.~Shamoyan, M.G.~Bashmakova}
\end{center}

\textbf{Abstract.}
{The intention of this survey to collect in one paper many recent results and advances related with Bergman type projection acting in various spaces of analytic functions in several complex variables in the unit ball, tubular domains over symmetric cones and bounded strongly pseudoconvex domains between function spaces of different dimensions.

Various new interesting extentions of old, classical results on Bergman projections will be provided in our survey.
Previously all these results were given in various papers of the first author. Bergman type projections have many nice applications in complex function theory of several complex variables in tubular domains over symmetric cones and in bounded strongly pseudoconvex domains. Our results can be seen as direct extensions of previously known results provided earlier by E. Stein, D. Bekolle, D. Debertol, B. F. Sehba, W. S. Cohn, C. Nana, L. Chen, Sh. Zhang and others which relate function spaces with the same dimension.

Practically all results on Bergman type projections acting between function spaces with different dimension may be valid both in context of tubular and pseudoconvex domains with similar proofs though we formulate some results only in tube or pseudoconvex domains.
Some our results are new even in case of simplest unit disk and polydisk. 

Problems which are related with Bergman type projection acting between various analytic and measurable function spaces of different dimensions in various complicated domains in $\mathbb{C}^n$ are as far as we know completely new and may have various new interesting applications in compex function theory. 

In this paper various new interesting problems in this research area will also be formulated and posed by authors. 

Our results may be valid also with very similar proofs in various Siegel domains of second type, matrix domains, bounded symmetric domains and various other domains in $\mathbb{C}^n$ with complicated structure.

Finally we remark for readers that we omit several new interesting results on Bergman type projections in analytic function spaces of several varables obtained by first author, hovewer we indicate all these papers in the list of references of this expository paper. }\\

\overfullrule=0pt \textbf{Keywords:} {Analytic function spaces, BMOA type spaces, Bergman ball, Bergman kernel, Bergman-type projections, Bergman-type spaces, Hardy spaces, Hardy-type spaces, Herz-type spaces, Trace theorems, Tubular and pseudoconvex domains, analytic functions, mixed norm spaces, product domains, unit ball. }\\

\section{Introduction, basic facts and assertions on tubular and bounded pseudoconvex domains}

Bergman projection is a classical topic in complex function theory of one and several complex variables.
Various aspects and applications of this research area are well known in complex function theory of one and several complex variables
New interesting applications of these type integral operators were found recently by first author.
In particular direct extentions of these integral Bergman type operators were defined in bounded strongly pseudoconvex domains and tubular domains over symmetric cones and they were applied in various trace theorems in particular.
These applications also concern various results connected with the extremal problems and trace theorems in various spaces of analytic functions of one and several variables and various complicated domains. This survey collect all these results from various papers of the first author in one paper.
We also discuss some new interesting open problems directly related with this topic.

In addition in his recent two decades the first author defined various new interesting scales of function spaces in product domains, spaces of several complex variables with the mixed norm in the unit ball, tubular domain and bounded strongly pseudoconvex domains, and new analytic Herz spaces in such type domains and studied the action of Bergman projections or the action of it is direct generalizations
in such type analytic function spaces of several complex variables.
These results may have separate interest for experts working in the area of complex function spaces and may be applied also for solutions of various interesting open problems of complex function theory in the future.
We also collect such type new interesting results on these integral operators obtained recently by the first author in this paper.

Various interesting new results on Bergman projections and Bergman type projections or it is direct generalizations were obtained recently also in very general Siegel domains of second type by various authors in various analytic spaces of several variables of Bergman type.
We refer the reader for this results to
[24] and [26] and various references there.
Some results are known for the action of Berman projection in multyfunctional spaces or multyfunctiinal expressions.
We refer the reader for example to [12] and [13, 43] and also various references there.
These all results were applied later to get so called new sharp decomposition theorems for multyfunctional
analytic function spaces in several complex variables.
We refer the reader to [39], [22], [23], [30], [31], [32] and various referenes there are also for this interesting valuable topic of complex function theory related with the Bergman projection. \\

For classical definitions of objects discussed below we referthe reader the classical textbooks of function theory in polydisk and unit ball
We provide some old, well known facts.
$T_{\beta}$ Bergman projection is bounded operator for large enough beta from $L^p_{\alpha}$ to $A^p_{\alpha}$ Bergman spaces, for all $p>1$, $\alpha>-1$
In the simple domains as unit ball and in the unit disk and in the unit polydisk.
Many interesting nice extensions of these not difficult results to analytic Bloch type spaces, to BMOA type spaces and Hardy type spaces are also well known to experts in context of mentioned rather simple domains. Later in recent two decades.
In connection with so called diagonal map direct extentions of many of these results to pair of spaces in different dimensions are also known in the polydisk and polyball.
It is a natural interesting new question to move further to extend these already known results to difficult domains such as tubular domains over symmetric cones and bounded strongly pseudoconvex domains with smooth boundary and to pair of function spaces of several variables defined in these very complicated domains and at the same time simultanuosly to function spaces different dimension.
This is exactly the main topic of this paper.

We define below basic objects and provide basic facts of function theory in tubular domains over symmeyric cones below and then in bounded strongly pseudoconvex domains which are very important for this paper and are needed for formulation of all our main theorems on projections in tubular and bounded strongly pseudoconvex domains.

A subset $\Omega$ of $\mathbb{R}^n$ or $V$ so that $dim V = n$ to be a cone if $\lambda x\in \Omega,$ for all $x\in\Omega, \lambda>0,$ if $\lambda x+\mu y\in\Omega$ for all $x,y\in\Omega, \lambda,\mu >0$ then it is convex. Let in
addition $\Omega^*=\{y\in \mathbb{R}^n : (y /x)>0$ for all $x\in \Omega\backslash \{0\}\}$ and $\Omega^*=\Omega.$ This type open
cone is selfdual ($\Omega^*$ is dual cone).\\
Let $G(\Omega)=\{g\in Gl(\mathbb{R}^n): g\Omega=\Omega\},$ where $Gl(\mathbb{R}^n)$ enotes the group of all linear
invertible transformation of $\mathbb{R}^n.$ If for all $x,y\in \Omega y=g x,$ for some $g\in G(\Omega)$ then
our open convex cone $\Omega$ is homogeneous, if also $\Omega=\Omega^*$
then it is symmetric cone.\\
Let $dv(w)$ and $dv_\alpha(w)=[\Delta^{\alpha-\frac{n}{r}}(v)]dudv, \alpha>\frac{n}{r}-1$ be a standard Lebesque measure in tubular
domains over symmetric cone $T_\Omega$ and weighted Lebesque measure in tube, $w=u+iv.$\\

Let $T_\Omega=V+i\Omega$ be the tube domain over an irreducible symmetric cone $\Omega$ in
the complexification $V^{\mathbb{C}}$
 of an $n$-dimensional Euclidean space $V$. $\mathcal{H}(T_\Omega)$ denotes
the space of all holomorphic functions on $T_\Omega$. We denote the rank of the cone $\Omega$ by $r$ and by $\Delta$ the determinant function on $V$.

For $\tau\in \mathbb{R}_+$ and the associated determinant function $\Delta(x)$ we set
$$A_{\tau}^{\infty}(T_\Omega)=\{F\in \mathcal{H}(T_\Omega): \|F\|_{A_{\tau}^{\infty}}=\sup_{x+iy\in T_\Omega}|F(x+iy)|\Delta^{\tau}(y)<\infty \} $$
It can be checked that this is a Banach space. For $1\leq p,q<\infty$ and $\nu\in\mathbb{R}, \nu>-1$ we denote by $A_\nu^{p,q}(T_\Omega)$ the mixed-norm weighted Bergman space consisting of analytic functions $f$ in $T_\Omega$ such that
$$\|F\|_{A_\nu^{p,q}}=\bigg(\int\limits_\Omega\big( \int\limits_{V}|F(x+iy)|^p dx\big)^{\frac{q}{p}}\Delta^\nu(y)dy \bigg)^{\frac{1}{q}}<\infty.$$
This is a Banach space. Replacing above $A$ by $L$ we will get as usual the corresponding larger space of all measurable functions in tube over symmetric cone with the same quazinorm. It is known that the $A_{\nu}^{p,q}(T_\Omega)$ space is nontrivial if
and only if $\nu>-1.$ When $p=q$ we write $A_\nu^{p,q}(T_\Omega)=A_\nu^p(T_\Omega).$ This is the classical weighted Bergman space with usual modification when $p=\infty.$

The (weighted) Bergman projection $P_\nu$ is the orthogonal projection from the
Hilbert space $L_\nu^2(T_\Omega)$ onto its closed subspace $A_\nu^2(T_\Omega)$ and it is given by the following
integral formula
$$ P_\nu f(z)=C_\nu\int\limits_{T_\Omega}B_\nu(z,w)z(w)dV_\nu ,$$
where
$$B_\nu(z,w)=C_\nu \Delta^{\nu+\frac{n}{r}}((z-\overline{w})/i)$$
is the Bergman reproducing kernel for $A_\nu^2(T_\Omega).$

Here we used the notation
$$dV_\nu(w)=\Delta^{\nu-\frac{n}{r}}(v)dudv.$$

Below and here we use constantly the following notations $w=u+iv\in T_\Omega$ and also $z=x+iy\in T_\Omega.$ Hence for any analytic function from $A_\nu^2(T_\Omega)$ the following integral formula is valid
$$ f(z)=C_\nu \int\limits_{T_\Omega} B_v(z,w)f(w)d_\nu w.$$

Below first we give basic definitions and define new and known analytic function spaces and standard objects in tubular domains over symmetric cones then provide such type definitions but incontext of bounded strongly pseudoconvex domain with smooth boundary.

Let $T_\Omega \in \mathbb{C}^n$ be a bounded tubular domains over symmetric cones in $\mathbb{C}^n$
We shall use the following notations:

$\bullet$ $\delta:T_\Omega \rightarrow \mathbb{R}^+$ will denote the determinant function from the boundary,
that is $\delta(z)=\Delta(Im z).$ Let $ d\nu_t(Z)=(\delta(z))^td\nu(z), t>-1;$

$\bullet$ $\nu$ will be the Lebesgue measure on $T_\Omega;$

$\bullet$ $H(T_\Omega)$ will denote the space of holomorphic function on $T_\Omega$ endowed
with the topology of uniform convergence on compact subsets;

$\bullet$ $B: T_\Omega\times T_\Omega \rightarrow \mathbb{C} $will be the Bergman kernel of $T_\Omega.$. Note that if $B$ is
kernel of type $t, t\in \mathbb{N},$ then $B^s$ is kernel of type $st, s\in\mathbb{N}, t\in\mathbb{N}.$

$\bullet$ Given $r\in(0, \infty)$ and $z_0\in T_\Omega$ we shall denote by $B_{T_\Omega}(z_0,r)$ the Bergman ball.

We define new Banach mixed norm analytic Bergman-type spaces in $T_\Omega\times...T_\Omega$ in product of tubular domains over symmetric cones as follows. Let $m\geq 1, p_j\in (1,\infty); \nu_j>\frac{n}{r}-1, j=1,...m.$
$$A_{\vec{\nu}}^{\vec{p}}=\{f\in H(T_\Omega)^m=H(T_\Omega\times ... \times T_\Omega) : $$
$$ \bigg( \int\limits_{T_\Omega}...\bigg( \int\limits_{T_\Omega} |f(z_1,...,z_m)|^{p_1}\Delta^{\nu_1-\frac{n}{r}}(y_1) dx_1dy_1\bigg)^{\frac{p_2}{p_1}}...\Delta^{\nu_m-\frac{n}{r}(y_m)dx_mdy_m} \bigg)^{\frac{1}{p_m}}<\infty )\}.$$

Let $\mathbb{C}$ denote the set
of complex numbers and let $\mathbb{C}^n=\mathbb{C}\times ... \times\mathbb{C}$
 denote the Euclidean space of complex dimension $n$. The open unit ball in $\mathbb{C}$ is the set $B_n=\{z\in \mathbb{C}: |z|<1\}$. We denote
by $H(B_n)$ the space of holomorphic functions on the open unit ball in $\mathbb{C}^n.$ Moreover,
let $\nu$ denote the Lebesgue measure on $B_n$ normalized in the sense that $\nu(B_n)=1$ and
for any $\alpha\in\mathbb{R},$ let $d\nu_{\alpha}(z)=c_\alpha(1-|z|^2)^\alpha\nu(z)$ for $z\in B_n.$ Here if $\alpha\leq -1, c_\alpha=1$ and if
 $\alpha>-1, c_\alpha=\frac{\Gamma(n\Gamma+\alpha+1)}{\Gamma(n+1)\Gamma(\alpha+1)}$ is the normalizing constant so that $\nu_\alpha$ has unit total
mass. The Bergman metric on $B_n$ is $\beta(z,w)=\frac{1}{2}\frac{1+|\varphi_z(w)|}{1-|\varphi_z(w)|},$ where $\varphi_z $ is the M\"{o}bius
transformation of $B_n$ that interchanges $0$ and $z.$ Let $D(a,r)=\{z\in B_n: \beta(z,a)<r \}$ denote the Bergman metric ball centered at $a\in B_n$ with radius $r > 0.$\\
We shall use the following notations:\\
$\bullet$ $\delta:D\rightarrow \mathbb{R}^+$ will denote the Euclidean distance from the boundary, that is $\delta(z)=d(z,\partial D);$\\
$\bullet$ given two non-negative functions $f,g: D\rightarrow \mathbb{R}^+$ we shall write $f\preceq g$ to say
that there is $C > 0$ such that $f(z)\leq Cg(z)$ for all $z\in D.$ The constant $C$ is
independent of $z\in D$, but it might depend on other parameters $(r,\theta, etc);$\\
$\bullet$ given two strictly positive functions $f,g: D\rightarrow \mathbb{R}^+$ we shall write $f\approx g$ if $f\preceq g$ and $g\preceq f,$
that is if there is $C > 0$ such that $C^{-1}g(z)\leq f(z)\leq Cg(z)$ for all $z\in D;$\\
$\bullet$ $\nu$ will be the Lebesque measure;\\
$\bullet$ $H(D)$ will denote the space of holomorphic functions on $D$, endowed with the
topology of uniform convergence on compact subsets;\\
$\bullet$ given $1\leq p\leq \infty,$ the Bergman space $A^p(D)$ is the Banach space $L^p(D)\cap H(D),$ endowed with the $L^p$
-norm;\\
$\bullet$ more generally, given $\beta\in\mathbb{R}$ we introduce the weighted Bergman space
$$A^p(D,\beta)= L^p(\delta^\beta \nu)\cap H(D)$$
endowed with the norm
$$\|f\|_{p,\beta}=\big [\int\limits_D |f(\zeta)|^p\delta^\beta(\zeta) d\nu(\zeta)\big]^{\frac{1}{p}}$$
if $1 \leq p <\infty,$ and with the norm
$$\|f\|_{\infty,\beta}=\|f \delta^\beta\|_{\infty}$$
if $p=\infty;$\\
$\bullet$ $K:D\times D\rightarrow\mathbb{C} $ will be the Bergman kernel of $D;$ the $K_t$ is a kernel of type $t;$\\
$\bullet$ for each $z_0\in D$ we shall denote by $k_{z_0}:D\rightarrow \mathbb{C}$ the normalized Bergman
kernel defined by
$$ k_{z_0}(z)=\frac{K(z,z_0)}{\sqrt{K(z_0,z0)}}=\frac{K(z,z_0)}{\|K(\cdot ,z_0)\|_2};$$
$\bullet$ given $r\in (0,1)$ and $z_0\in D,$ we shall denote by $B_D(z_0,r)$ the Kobayashi ball of
center $z_0$ and radius $\frac{1}{2}\log \frac{1+r}{1-r}.$

For any two $m$-tuples of real numbers $a=(a_1,...,a_m)$ and $b=(b_1,...,b_m),$ we define the integral operator in the unit ball and products of balls.
$$(S_{a,b}f)(z_1,...,z_m)=\prod\limits_{j=1}^m (1-|z_j|^2)^{a_j}\int
\limits_{B_n} \frac{f(w)(1-|w|^2)^{-n-1+\sum_{j=1}^m b_j}}{\prod_{j=1}^m (1-\langle z_j,w\rangle)^{a_j+b_j}}d\nu(w),$$
where $z_1,...z_m$ are in $B_n,$ and $f$ is a function in $L^1(B_n, d\nu_{-n-1-\sum_{j=1}^m b_j}).$ Note that for
such $f$ the function $S_{a,b}f$ is defined on $(B_n)^m,$ the product of $m$ copies of $B_n$.
These new extensions of Bergman projections $S_{a,b}$ Bergman type integral operators where used and studied in trace type theorems in BMOA type analytic function spaces by first author in the ball tube and bounded strongly pseudoconvex domains.
We can easily based on definition above define such type $S_a$,b integral operators in context of all mentioned domains in $\mathbb{C}^n$.

% Свойство C.
For our theorem in pseudoconvex domains ($p < 1$ case) for
weighted Bergman kernels we always assume a little bit more, namely the following
properity (we call it (C)) $K_t(z,w)\asymp K_t(a_k,w)$ for any $z\in B_D(a_k,r), r\in(0;1), w\in B_D(a_m,r), r\in(0;1)$ for $t>0,$ and any natural number $m$(we assume here the Bergman kernel is positive, otherwise add modulus.)

Let $T_\Omega$ be the tube domain over symmetric cone, and $H(T_\Omega)$ be the space of all analytic functions in tube. We define
Bergman spaces for $1\leq p,q<\infty, \gamma>\frac{n}{r}-1.$ Let
$$A_{\gamma}^{p,q}(T_\Omega)=$$
$$=\{f\in H(T_\Omega) : \bigg( \int\limits_\Omega \bigg( \int\limits_{\mathbb{R}^n} |f(x+iy)|^p dx \bigg)^{\frac{q}{p}} (\Delta^{\gamma-\frac{n}{r}}(y)dy) \bigg)^{\frac{1}{q}}<\infty\}.$$

Replacing $A$ by $L$ and $H$ by $L^1$ we get as usual known larger spaces of
measurable functions in tube $T_\Omega.$ Note $A_\gamma^{p,q}=\{0\}$ if $\gamma<\frac{n}{r}-1.$

%We would like to remark here that $A^{p,q}_\alpha$ spaces of analytic functions which defined above can be extended by adding the amount of integrals rather easily in the norm of a space in two different ways to product domains and these are new function spaces and new various projection theorems in these function spaces will be interesting also to study and to obtain. We leave this interesting problems to readers.
%Дублирует нижнее, перед \section

Let $f\in H(T_\Omega\times...\times T_\Omega), f=(z_1,...,z_m), m\geq 1;$ this means f function is holomorphic separately
by each variable $(z_i).$ Let $\vec{f}=(f_1,...,f_m); f=\prod\limits_{j=1}^m f_j,$
$$ \|\vec{f}\|^{p}_{B_{\vec{\nu}}^p}=\int\limits_{T_\Omega} |f_1(z)|^p...|f_m(z)|^p\Delta^{\nu_1-\frac{n}{r}}(y)...\Delta^{\nu_m-\frac{n}{r}}(y)dxdy<\infty;$$
$y=Im z; f_j\in H(T_\Omega), 1\leq p<\infty, \nu_j>\frac{n}{r}-1,j=1,...,m. $ We assume $ \|\vec{f}\|^{p}_{B_{\vec{\nu}}^p}<\infty.$

These are multifunctional analytic function spaces in tube. We can define similarly them in bounded strongly pseudoconvex domains. Many theorems of this papers can be passed to such type spaces also by similar methods.

Next let $1\leq p<\infty, f=f(z_1,...,z_m),$ we consider analytic subspaces
of $H(T_\Omega^m), T_\Omega^m=T_\Omega\times...\times T_\Omega, \nu_j>\frac{n}{r}-1, j=1,...,m.$ These are spaces $(A_\nu^p)_1, (A_\nu^p)_2, (A_\nu^p)_3$ with norms
$$\|f\|^p_{(A_{\vec{\nu}}^p)_1}=\int\limits_{T_\Omega}...\int\limits_{T_\Omega} |f(x_1+iy_1,...,x_m+iy_m)|^p \prod\limits_{j=1}^m \Delta^{\nu_j-\frac{n}{r}}(y_j)dx_j dy_j<\infty;$$
$$\|f\|^p_{(A_{\vec{\nu}}^p)_2}=\int\limits_{\mathbb{R}^n}...\int\limits_{\mathbb{R}^n}\int\limits_\Omega |f(x_1+iy_1,...,x_m+iy_m)|^p \Delta^{\nu-\frac{n}{r}}(y) \bigg(\prod \limits_{j=1}^m dx_j \bigg)dy< \infty;$$
$$\|f\|^p_{(A_{\vec{\nu}}^p)_3}=\int\limits_{\mathbb{R}^n}\int\limits_{\Omega}...\int\limits_\Omega |f(x_1+iy_1,...,x_m+iy_m)|^p \prod_{j=1}^m \Delta^{\nu_j-\frac{n}{r}}(y_j) dxdy_j<\infty.$$

%Кусок текста из промежутка между 19 и 20 теоремой
Let $D$ be a bounded strongly pseudoconvex domain in $\mathbb{C}^{n}$
with smooth boundary, let $d(z)=dist(z, \partial D).$ Then there is a neighbourhood $U$ of $\bar{D}$ and $\rho \in \mathbb{C}^{\infty}(U)$ such that $D=\{z\in U : \rho(z)>0\}, |\Delta \rho(z)|\geq c>0$ for $z\in\partial D, 0<\rho(z)<1$ for $z\in D$ and $-\rho$ is strictly plurisubharmonic in a neighborhood $U_0$ of $\partial D.$

Then there is an $r_0$ such that the domains $D_r=\{z\in D: \rho(z)> r\}$ are also smoothly
bounded strictly pseudoconvex domains for all $0<r\leq r_0.$ Let $d\sigma_r$ be the normalized surface
measure on $\mathbb{C}^{n}$. For $0<p<\infty, 0<q\leq \infty, \delta>0, k=0,1,r$ we set
$$\|f\|_{p,q,\delta,k}=\sum\limits_{|\alpha|\leq k}\bigg( \int\limits_0^{r_0}\bigg(r^\delta \int\limits_{\partial D_r}|D^\alpha f|^p d\sigma_r\bigg)^{\frac{q}{p}} \frac{dr}{r} \bigg)^{\frac{1}{q}}, 0<q<\infty;$$
$$\|f\|_{p,\infty,\delta,k}=\sup\limits_{0<r< r_0}\sum\limits_{|\alpha|\leq k} \bigg(r^\delta \int\limits_{\partial D_r}|D^\alpha f|^p d\sigma_r\bigg)^{\frac{1}{p}}.$$
The corresponding $ A_{\delta,k}^{p,q}=\{f\in H(D): \|f\|_{p,q,\delta,k} < \infty\}$ are complete quasinormed spaces
for $p,q\geq 1.$ They are Banach spaces. We mostly deal with the case $k=0$ then we write
simply $A_\delta^{p,q}$ and $\|f\|_{p,q,\delta}.$

We also consider these spaces for $p=\infty$ and $k=0,$ the corresponding space is denote by $A_\delta^{\infty,p}=A_\delta^{\infty,p}(D)$and consists of all $f\in H(D)$ such that
$$\|f\|_{\infty,p,\delta}=\bigg(\int\limits_0^{r_0} \big( \sup\limits_{\partial D_r}|f|\big)^p r^{\delta p-1}dr\bigg)^{\frac{1}{p}}<\infty.$$
Also for $\delta>-1$ the space $A_\delta^\infty =A_\delta^{\infty}(D)$ consists of all $f\in H(D)$ such that
$$ \|f\|_{A_\delta^{\infty}}=(\sup\limits_{z\in D})|f(z)|\rho^{\delta}(z)<\infty,$$
and the weighted Bergman space $A_\delta^p=A_{\delta+1}^{p,p}(D)$ consists of all $f\in H(D)$
such that
$$ \|f\|_{A_\delta^p}=\bigg( \int\limits_{D}|f(z)|^p \rho ^{\delta}(z)d\nu(z)\bigg)^{\frac{1}{p}}<\infty,$$
where $d\nu$ is a normalized Lebegues measure in $D$.

%Кусок текста из промежутка между теоремами 20 и 21
Let $T_\Omega=v+i\Omega$ be the tube domain over an irreducible symmetric cone $\Omega$ in the
complexification $V^{\mathbb{C}}$ of an $n$-dimensional Euclidean space $\tilde{V}.$ We denote the rank of the cone $\Omega$ by $r$ and by $\Delta$ the determinant function on $\tilde{V}.$

Let us introduce some convenient notations regarding multi-indices.
If $t=(t_1,...,t_m)$ then $t^*=(t_m,...,t_1)$ and, for $a\in\mathbb{R}, t+a=(t_1+a,...,t_m+a).$ Also, if $t,k\in \mathbb{R}^m,$ then $t<k$ means $t_j<k_j,$ for all $j=1,..,m.$\\
We are going to use the following multi-index
$$ g_0=\big( (j-1)\frac{d}{2}\big)_{1\leq j\leq r}, where (r-1)\frac{d}{2}=\frac{n}{r}-1.$$
$\mathcal{H}(T_\Omega)$ denotes the space of all holomorphic functions on $T_\Omega$. We denote $m$ Cartesian
products of tubes by $T_\Omega^m=T_\Omega\times...\times T_\Omega$ the space of all analytic function on this
new product domain which are analytic by each variable separately will be denoted by $\mathcal{H}(T_\Omega^m).$ By $m$ here and everywhere below we denote below a natural number bigger
than 1. For $\tau\in \mathbb{R}_+$ and the associated determinant function $\Delta(x)$ we set
$$ \tilde{A}_\tau^{\infty}(T_{\Omega}) = \bigg \{ F\in \mathcal{H}(T_\Omega) : \|F\|_{A_\tau^{\infty}}=\sup\limits_{x+iy\in T_\Omega} |F(x+iy)|\Delta^\tau (y)<\infty \bigg \}.$$
It can be checked that this is a Banach space.\\
For $1\leq p,q <+\infty$ and $\nu\in \mathbb{R}, \nu>\frac{n}{r}-1$ we denote by $\tilde{A}_{\nu}^{p,q}(T_\Omega)$ the mixed-norm
weighted Bergman space consisting of analytic functions $f$ in $T_\Omega$ that
$$ \|F\|_{\tilde{A}_\nu^{p,q}}= \bigg ( \int\limits_\Omega \bigg( \int\limits_{\tilde{V}} |F(x+iy)|^p dx \bigg)^{\frac{q}{p}} \Delta^{\nu} (y)\frac{dy}{\Delta(y)^{\frac{n}{r}}} \bigg)^{\frac{1}{q}}<\infty.$$
This is a Banach space.\\
Replacing above simply $\tilde{A}$ by $\tilde{L}$ we will get as usual the corresponding larger space
of all measurable functions in tube over symmetric cone with the same quazinorm. It is known the $\tilde{A}_{\nu}^{p,q}(T_\Omega)$ space is nontrivial if and only if $\nu>\frac{n}{r}-1.$ And we will assume this everywhere below. When $p=q$ we write $\tilde{A}_{\nu}^{p,q}(T_\Omega)=\tilde{A}_{\nu}^{p}(T_\Omega).$ This is the classical weighted Bergman space with usual modification when $p=\infty$.
We add some notions on Bergman type analytic function spaces on products of tubular
domains.

Let $T_\Omega^m=T_\Omega\times...\times T_\Omega.$ To define related two Bergman-type spaces $\tilde{A}_\nu^p(T_\Omega)$ and $\tilde{A}_\tau^{\infty}(T_\Omega)$ ( $\nu$ and $\tau$ can be also vectors) in $m$-products of tube domains $T_\Omega^m$ we follow
standard procedure which is well-known in case of unit disk and unit ball. Namely we consider analytic $F$ functions $F=F(z_1,...,z_m)$ which are
analytic by each $z_j, j=1,...,m$ variable, and where each such variable belongs to $T_\Omega$ tube and define as $H(T_\Omega^m)$ the space of all such functions. For example we set, for all $z_j=x_j+i y_j, \tau_j\in \mathbb{R}, j=1,...,m, F(z)=F(z_1,...,z_m), \tau=(\tau_1,...,\tau_m)$
$$\tilde{A}_\tau^{\infty}(T_\Omega^m)=\bigg \{ F\in\mathcal{H}(T_\Omega): \|F\|_{\tilde{A}_\tau^{\infty}}= \sup\limits_{x+iy\in T_\Omega^m}|F(x+iy)|\Delta^\tau(y)<\infty \bigg\},$$ where
$$|F(x+iy)|=|F(x_1+iy_1,...,x_m+iy_m)|,$$ and $\Delta^\tau(y)$ is a product of $m$ one-dimensional $\Delta^{\tau_j}(y_j)$ functions, $j=1,...,m.$
Similarly the Bergman space $\tilde{A}_\tau^p$ can be defined on products of tubes for all $\tau=(\tau_1,...,\tau_m), \tau_j>\frac{n}{r}-1, j=1,...,m.$ It can be shown that all spaces are Banach spaces. Replacing above simply $\tilde{A}$ by $\tilde{L}$ we will get as usual the corresponding larger space of all measurable functions in products of tubes over symmetric cone with the same quazinorm.

The (weighted) Bergman projection $P_\nu$ is the orthogonal projection from the Hilbert
space $\tilde{L}_\nu^2(T_\Omega)$ onto its closed subspace $\tilde{A}_\nu^2(T_\Omega)$ and it is given by the following integral
formula (see (\cite{d2}))
$$P_\nu f(z)=C_\nu \int\limits_{T_\Omega}B_\nu(z,w) f(w)d\tilde{V}_\nu(w),$$ where
$$B_\nu(z,w)=C_\nu \Delta^{-\nu+\frac{n}{r}}(\frac{z-\bar{w}}{i})$$
is the Bergman reproducing kernel for $\tilde{A}_\nu^2(T_\Omega).$\\
Here we used the notation $d\tilde{V}_\nu(w)=\Delta^{\nu-\frac{n}{r}}(v)du dv.$ Below and here we use constantly the following notations $w=u+iv \in T_\Omega$ and also $z=x+iy \in T_\Omega.$\\
Hence for any analytic function from $\tilde{A}_\nu^2(T_\Omega)$ the following reproducing integral
formula is valid
$$ f(z)=C_\nu\int\limits_{T_\Omega} B_\nu(z,w) f(w) d\tilde{V}_\nu(w).$$
We provide now several results in this area on Bergman type projections in tube.
The weighted Bergman projection $P_\nu$ is the orthogonal projection from the Hilbert
space $\tilde{L}_\nu^2(T_\Omega)$ onto its closed subspace $\tilde{A}_\nu^2(T_\Omega)$ and it is given by the integral formula
$$ (P_\nu f)(z)=\int\limits_{T_\Omega}B_\nu(z,w) f(w)\Delta^{\nu-\frac{n}{r}}(Im w)dv(w),$$
$z\in T_\Omega$ and $\nu>\frac{n}{r}.$\\
If $P_\nu$ extends to a bounded operator on $\tilde{L}_\nu^{p,q}$ then the topological dual space $(\tilde{A}_\nu^{p,q})^*$ of the Bergman
space $tilde{A}_\nu^{p,q}$ identifies with $tilde{A}_\nu^{p',q'}$ under the integral pairing
$$ \langle f,g \rangle_\nu =\int\limits_{T_\Omega} f(z)\overline{g(z)} \Delta^{\nu-\frac{n}{r}}(Im z)dv(z),$$
$f\in \tilde{A}_\nu^{p,q}, g\in \tilde{A}_\nu^{p',q'} $ (see (\cite{d4})). Let $\beta>-1, \gamma>0$ and $\alpha>0.$ Let also
$$(T_{\alpha,\beta,\gamma}f)(z)=\Delta^{\alpha}(Im z) \int\limits_{T_\Omega} B_\gamma(z,w) f(w) \Delta^{\beta}(Im w) dv(w),$$
$$(T^+_{\alpha,\beta,\gamma}f)(z)=\Delta^{\alpha}(Im z) \int\limits_{T_\Omega}| B_\gamma(z,w) |f(w) \Delta^{\beta}(Im w) dv(w),$$
$z\in T_\Omega, f\in \tilde{L}^1(T_\Omega).$

We also define new Herz type integral operator for positive $\alpha_j, j=1,...,m$ and for all $\beta>-1, \gamma>-1$
$$[T_{\vec{\alpha},\beta,\gamma}(g)](z_1,...,z_m)=\int\limits_{T_\Omega}\int\limits_{B(\tilde{w},R)}\frac{g(w)[\Delta^\beta(Im w)]dv(w)}{\big[\prod_{j=1}^m \Delta^{\alpha_j}(\frac{z_j-\bar{w}}{i})\big]}\Delta^{\gamma} (Im \tilde{w})dv(\tilde{w}),$$
where $z_j\in T_\Omega, j=1,...,m.$
%Здесь куски текста кончаются

These new objects, quazinorms and various function spaces which we defined in this section in tubular domains over symmetric cones can be also easily similarly defined in he unit ball and in bounded strongly pseudoconvex domains with the smooth boundary.
This can be done by interested readers. We would like to remark here that $A_{\alpha}^{p,q}$ spaces of analytic functions which defined above can be extended by adding the amount of integrals rather easily in the norm of a space in two different ways to product domains and these are new function spaces and new various projection theorems in these function spaces will be interesting also to study and to obtain. We leave this interesting problems to readers.

Many Bergman projection theorems acting between spaces of different dimension appeared in papers of the first author as a part of solution of various trace problems in function spaces previously.
We also refer the reader for such results to papers of the first author related with trace problem posted in various databases.

%Возможно куски текста стоит поместить сюда

We alert the reader some objects which can be seen in our theorems are not defined in this note, we refer the reader for them to related and mentioned in those theorems research papers.

In this paper various new interesting problems in this research area will also bе formulated and posed by authors.

Complete analogues of all function spaces we defined in this section in the tubular domains can be defined in pseudoconvex domains easily by readers and the reverse is also true.

Practically all our results were known previously in very simple domains and for particular values of parameters.

We denote as usual in this paper by $C$ with various lower indexes is Bergman representation constant.

All theorems of this paper may also have with similar proofs complete analogues in Bloch-type spaces in tube and pseudoconvex domains which we defined above.

\section {Bergman type projections in tubular domains and bounded pseudoconvex domains}

The main intention of this large second section to collect many recent results mainly related with advances of first author and his coauthors concerning Bergman projection in various analytic function spaces in rather complicated tubular domains over symmetric cones and bounded stongly pseudoconvex domains with smooth boundary in $\mathbb{C}^n$.
We at the same time sometimes in one theorem generalize certain known classical results simultanuosly in two directions, which add some additional interest to such type assertions.
We pay attention to certain specific aspects of issues related with Bergman projection which often appear in relation with the so called trace problem in analytic function spaces of several variables $\mathbb{C}^n$.
Namely we are interested on Bergman type projections acting between spaces of measurable and analytic functions of different dimension. We give a very simple example below for readers first.

Let $X$ be a space of measurable functions on a fixed $D$ domain in $\mathbb{C}^{n}$. The question is there a Bergman type projection acting between this space and another space of analytic function but defined on $D \times ... \times D$ domain. Which at the same time generalize time an old well-known classical result.
Note actions of such type Bergman type integral operators from various interesting spaces of measurable functions defined on a certain fixed $D$ domain in $\mathbb{C}$ or $\mathbb{C}^n$ to such type analytic spaces defined also on $D$ domain namely such type theorems on Bergman type projection between various function spaces defined on the same $D$ domain are well known for many experts and can be seen in literature for simple domains as the unit disk, polydisk, and unit ball and even for various more difficult domaims in $\mathbb{C}^n$.
The interesting new issues here is to consider first difficult domains in $\mathbb{C}^{n}$. Then to extend these known already results as far as possible to $(X,Y)$ pairs of various interesting function spaces defined in $D$ and $D \times ... \times D$, namely $X(D)$ and $Y(D \times ... \times D)$.

We collect mainly such type new interesting Bergman projection theorems below. We hope that applications of such type new interesting results in analytic function spaces in tubular domains over symmetric cones and bounded strongly psuedoconvex domains with the smooth boundary may be numerous.

Below sometimes in our theorems we do not define related to these assertions function spaces describing them shortly before theorems and referying the interested reader to related paper indicated in formulation of those mentioned our theorems.

In two theorems below we give complete descriptions of traces of certain new function spaces in tubular domains over simmetric cones and within these results we provided new interesting results on the action on Bergman projection between function spaces of different dimensions in tube domains. 
We refer the reader for definitions in mentioned in these theorems Herz spaces in tube related papers which we mention there.

Bergman projections acting between function spaces with different dimension are usually closely related with so called trace problem in analytic function spaces of several complex variables. We say trace $X=Y$, when $X$ and $Y$ are quazinormed spaces of $H(D)$ and $H(D\times...\times D)$, for a certain $D$ domain in $\mathbb{C}^n$, if for any $f$, from $X$, $f(z,...,z)$ belong to $Y$ and the reverse is also true, any $g$ function from $Y$ can be extended to a fixed $f$ function, $f\in X$, so that $f(z,...,z) = g(z)$, $z\in D$.

\begin{theorem}\cite{0} % теорию из этой статьи я не написала
Let $f\in A_{\vec{\nu}^{\vec{p}}}(T_\Omega^m),1\leq p_j<\infty, \nu\in \mathbb{R}^m, \nu_j>\nu_0$ for fixed $\nu_0=\nu_0(n,p,r,m),$ for all $j=1,...,m.$ Then $f(z,...,z)\in A_s^{p_m}, s=\nu_m-\frac{n}{r}+\sum_{j=1}^{m-1}(\nu_j+\frac{n}{r})\frac{p_m}{p_j}$ with related estimates for norms. And for all $\tilde{p}_1, \frac{1}{\tilde{p}_1}+\frac{1}{p_m}=1, \frac{n}{r} \leq \tilde{p}_1, j=1,...,m,$ for some fixed $k$ large enough, $\nu_j>k$ the reverse is also true. For each function $g\in A_s^{p_m}(T_\Omega)$ there is an $F$ function $F(z,...,z)=g(z), F\in A_{\vec{\nu}}^{\vec{p}}(T_\Omega).$
Let in addition
$$T_\beta(f)(z_1,...,z_m)=C_\beta\int\limits_{T_\Omega} f(w) \big( \prod\limits_{j=1}^m \Delta^{-t}(\frac{z_j-\bar{w}}{i})\big)dV(w), mt=\beta+\frac{n}{r}, z_j\in T_\Omega, j=1,...m.$$
Then following assertion hold for all $\beta,$ so that $\beta>\beta_0$ for some fixed large enough
positive number $\beta_0.$ The $T_\beta$ Bergman-type integral operator (expanded Bergman projection) maps $A_s^{p_m}(T_\Omega)$ to $A_{\vec{\nu}}^{\vec{p}}(T_\Omega^m), \nu=(\nu_1,...,\nu_m), \nu_j>\nu_0, j=1,...,m.$
\end{theorem}

\begin{theorem}\cite{0}
Let $\nu_j>\nu_0$ and $\tau_j>\tau_0$ for some fixed positive numbers $\nu_0=\nu_0{p,n,r,m}$ and $\tau_0=\tau_0(p,n,r,m), 1\leq p<\infty. $ Let $f\in S^{p}_{\nu,\tau}(T_\Omega^m) $ then $f(z,...,z)$ belongs to $A_s^p(T_\Omega)$ where $s=\sum\limits_{j=1}^m (\nu_j+\frac{n}{r})p+\sum\limits_{j=1}^m (\tau_j+\frac{n}{r})-\frac {n}{r}$ and for every $f$ function $f\in A_s^p(T_\Omega)$ there is an $F$ function $F\in S_{\nu,\tau}^p$ so that $F(z,...,z)=f(z)$ for all $\frac{n}{r}\leq p', \frac{1}{p}+\frac{1}{p'}=1.$//
Let in addition
$$ (T_\beta f)(z_1,...,z_m)= C_\beta \int\limits_{T_\Omega} f(w) \prod\limits_{j=1}^m \Delta^{-t} \big(\frac{(z_j-\bar{w})}{i}\big)dV_\beta(w),$$
$mt=\beta+\frac{n}{r}, z_j\in T_\Omega, j=1,...,m.$ Then the following assertions holds for all $\beta$ so that $\beta>\beta_0$ for somefixed large enough positive number $\beta_0.$ The $T_\beta$ Bergman
projection-type integral operator maps $A_s^p(T_\Omega)$ to $S_{\nu,\tau}^p(T_\Omega^m), \nu=(\nu_1,...,\nu_m), \tau=(\tau_1,...,\tau_m), \nu_j>\nu_0, \tau_j>\tau_0, j=1,...,m.$
\end{theorem}

Proofs of theorems $1, 2$ have big amount of similarities with the much simpler and already known unit disk case. The role of so called $r$-lattices in tubular domains are very important in these proofs. These results with very similar proofs are probably valid in bounded strongly pseudoconvex domains with the smooth boundary, based on known properities of $r$-lattices in mentioned bounded strongly pseudoconvex domains.

Let $U$ be the unit disk, $T$ the unit circle.
We define as usual analytic and well studied by many experts  $A^{p,q}_\alpha$ spaces for all positive $p, q$, $\alpha>-1$, as function spaces with finite quazinorms.
$$\|f\|_{p,q,\alpha}=
\int\limits_0^1(\int\limits_T |f(r\xi)|^{p}d\xi)^{q/p}\times(1-r)^{\alpha}(dr).$$

Very similarly these analytic function spaces can be easily defined in polydisk 
Similarly these function spaces can be easily defined also in very general tubular and bounded pseudoconvex domains in $\mathbb{C}^n$ and even in products of such type domains.by  adding the amount of integrals in quazinorms.
We denote such spaces of measurable functions by $L^{p,q}_\alpha$ and by 
$L^{p,q}_\alpha$ putting wave on L dealing with spaces on product of such type domains.
And replacing $L$ by $A$ for analytic subspaces.
In the following theorems various projection results are given for such type function spaces.
The projection theorems between spaces of same dimension in tube and pseudoconvex domains and spaces on products of such type domains can be easily obtained by readers by rather simple repetition of arguments of proofs provided earlier in such type spaces in the unit disk case.
Other such type results between such spaces of different dimensions are more interesting, here proofs again are heavily based on arguments provided earlier by other authors in the polydisk.

%- Начинается 01 (без теории)
\begin{theorem}\cite{0-1}
Define the operator $$T_\beta f(z)=\int\limits_\Lambda f(w) K_{\beta+n+1}(z,w) (\delta^\beta(w))dv(w)$$
where $\beta$ is large enough. Then the operator $T_\beta$ maps the space $L_\delta^{p,q}$ into the analytic function space $A_\delta^{p,q}$ for all $1<p,q<\infty$ and $\delta>-1$.
\end{theorem}

\begin{theorem}\cite{0-1}
Consider the integral operator
$$ (\tilde{T}_\beta f)(z_1,...,z_m)= \int\limits_\Lambda ... \int\limits_\Lambda (\vec{f}(w))\bigg( \prod\limits_{j=1}^m K_{\beta+n+1}(z_j,w_j)\bigg)(\delta^\beta(w_j))dv(w_j).$$
Then for all for all $\beta>\beta_0$ with $\beta_0$ large enough, the operator $\tilde{T}_\beta$ maps the space $\tilde{L}_\delta^{p,q}$ into the space $\tilde{A}_\delta^{p,q}$ for all $z_j\in \Lambda, j=1,...,m, 1<p,q<\infty, $ and $\delta>0.$
\end{theorem}

\begin{theorem}\cite{0-1}
For all $1<p,q<\infty$ the operator
$$(S_\beta f)(z_1,...,z_m)=\int\limits_\Lambda (f(w))\bigg( \prod\limits_{j=1}^m K){\beta+n+1}(z_j,w)\bigg)((\delta^{\beta_1})(w))dv(w), $$
with $\beta>\beta_0$ and $\beta_0$ large enough, where $\beta_1=(\beta+n+1)m-(n+1),$ maps the space $(L_\delta^{p,q})(\Lambda)$ to $(\tilde{A}_\delta^{p,q})(\Lambda\times...\times\Lambda)$ if for any $F\in (\tilde{A}_{\delta_1}^{p',q'})$ the condition $ F(w,...,w)\in (A_{\delta_2}^{p',q'})$ is satisfied, where $ \delta_2=\delta(1-q')+\beta m q'+(n+1)(m-1)(q'), \delta_1=(\beta-\delta)q'+\delta. $

\end{theorem}

\begin{theorem}\cite{0-1}
Consider the integral operator
$$ (T_\tau f)(z)= \int\limits_{T_\Lambda}(f(\tilde{w})) B_{\tau+\frac{2n}{r}}(r,\tilde{w}) [\Delta^{\tau}(Im \tilde{w})]dv(\tilde{w}). $$
Then, for all $\tau>\tau_0$, with $\tau_0$ large enough, the operator $(T_\tau)$ maps $L_{\nu}^{p,q}(T_\Lambda)$ to $A_\nu^{p,q}(T_\Lambda)$, where
$$ (A_\nu^{p,q})=\bigg \{ f\in H(T_\Lambda) : \int\limits_\Lambda \bigg ( \int\limits_{R_n}|f(x+iy)|^p dx \bigg)^{\frac{q}{p}} \Delta^\nu (Im z)dxdy <\infty \bigg\},$$
where $\nu>-1$ and $1<p,q<\infty.$
\end{theorem}

\begin{theorem}\cite{0-1}
For every $1<p,q<\infty$, the operator
$$ (S_\beta f)(z_1,...,z_m)=\int\limits_\Lambda (f(w))\bigg(\prod_{j=1}^m B_{\beta+\frac{2n}{r}}(z_j,w) \bigg)((\delta^{\beta_1})(w))dv(w),$$ with $\beta>\beta_0,$ where $\beta_0$ large enough, and $\beta_1=(\beta+\frac{2n}{r})m-\frac{2n}{r},$ maps the space $(L^{p,q}_\alpha)(\Lambda)$ to $(\tilde{A}_\delta^{p,q})(\Lambda\times...\times\Lambda)$ if for any $F\in (\tilde{A}_{\delta_1}^{p',q'}),$ the condition $F(w,...,w)\in A_{\delta_2}^{p',q'}$ is satisfied, where
$$ \delta_2=\delta(1-q) +\beta mq'+(2n/r)(m-1)(q'), \delta_1=(\beta-\delta)q'+\delta.$$
\end{theorem}

For $m=1$ additional condition in theorem 5 and 7, on Bergman projection between spaces with different dimension vanishes and we get assertions which are well known in the unit ball and unit polydisk.
These new interesting results of theorem 5 and 7 may be valid also Siegel domains of second type, bounded symmetric domains and matrix domains and also in minimal bounded homogeneous domains with very similar proofs.
In theorem 8 below we formulate a new theorem on Bergman type integral operators which related spaces with different dimension in tubular domains and bounded strongly pseudoconvex domains, this theorem is well known for particular case, $m=1$ case in tube and pseudoconvex domains
The definition of $\mathbb{R}_x$ expanded Bergman projection can be seen in [7]
This theorem 8 also has complete analogues in various domains in $\mathbb{C}^n$ in particular bounded symmetric domains, Siegel domains of second type, matrix domains, etc with the same proof.

\begin{theorem}\cite{3} %6 n.19

1). Let $1<p<\infty, s_j>\frac{n}{r}-1, j=1,...,m.$
Then for some fixed large enough $X_0$ and all $X_j>X_0, j=1,...,m, $ there is a constant $C>0$ such that
$$\int\limits_{T_\Omega}|R_{\overrightarrow{X}}g(w)|^p \Delta^{(\sum\limits_{j=1}^{m})s_j+\frac{2n}{r}(m-1)}(Im w)dudv\leq \tilde{C}\| g\|_{L_{\overrightarrow{s}}^p}(T_\Omega^m).$$
2). Let $1<p<\infty, s_j>\frac{n}{r}-1, j=1,...,m. $

Then for some fixed large enough $X_0,$ and all $X_j>X_j, j=1,...,m$ there is a constant $C>0$
such that
$$\int\limits_{\Lambda}|\tilde{R}_{\overrightarrow{X}}g(w)|^p \delta^{\tau}(w)dv(w)\leq \tilde{C}\| g\|_{L_{\overrightarrow{s}}^p}(\Lambda^m),$$
here $\tau=\big(\sum\limits_{j=1}^m s_j\big)+(m-1)(n+1).$

\end{theorem}

Let further $dv_{\alpha}=(\delta^{\alpha})dv(z)$ where $\delta(z)=\rho(z).$ Since $|f(x)|^p$ is subharmonic (even plurisubharmonic) for a holomorphic $f$ we have $A_s^p(D)\subset A_t^{\infty}(D),$ for $0<p<\infty, sp>n$ and $t=s$ also $A_s^p(D)\subset A_s^1(D)$ for $0<p\leq 1$ and $A_s^p(D)\subset A_t^1(D)$ for $p>1$ and $t$ sufficiently large.

Therefore we have an integral representation
$$ f(z)=\int\limits_{D}f(\xi)K_{\tilde{t}}(z,\xi)\rho^t(\xi) dv(\xi),$$
for $f\in A_t^1(D), z\in D, \tilde{t}=t+n+1,$ where $K_{\tilde{t}}(z,\xi)$ is a kernel of type $t$ that is a
smooth function on $D\times D$ such that
$$ |K_{\tilde{t}}(z,\xi)|\leq C_1|\widetilde{\Phi}(z,\xi)|^{-(n+1+t)}, $$ where $\widetilde{\Phi}(z,\xi)$ is so called Henkin-Ramirez function on $D$.

In two theorems below formulated in context of bounded strongly pseudoconvex $D$ domains ($dv_\beta$ is a usual weighted Lebeques measure in $D$) we formulate direct extentions of known results in the polydisk to these general $D$ domains in $\mathbb{C}^n$.
These results with the same proof are valid in tubular domains and Siegel domain of second type and bounded symmetric domains and also in various matrix domains.

\begin{theorem}\cite{4} %7 n.20
Let $0<p,q\leq 1,\alpha>\tilde{\alpha}_0, \tilde{\alpha}_0 $ is large enough then
$$ (T_\alpha f)(w)= \int\limits_{\Lambda}K_{\alpha_0}(z,w)f(z)dv_{\alpha}(z), \alpha_0=\alpha+n+1,$$
maps $A_\beta^{p,q}$ into $A_\beta^{p,q}$ for all $\beta, \beta>0, dv_\alpha=\delta^{\alpha}dv.$
\end{theorem}

\begin{theorem}\cite {4} %8 n.21
Let $p_i>1, i=1,...,m,\alpha_j>-1, j=1,...,m.$. Then we
have that the $V_{\overrightarrow{\beta}}(f)$ operator, where
$$V_{\overrightarrow{\beta}}(f)=\int\limits_D ... \int\limits_D f(z_1,...,z_m) \prod\limits_{j=1}^m K_{\beta_j+n+1}(z_j,w_j) dv_{\beta_1}(z_1) ... dv_{\beta_m}(z_m),$$
for all $\beta_j>\beta_0, j=1,...,m, \beta_0$ is large enough, maps $L_{\alpha_1,...,\alpha_n}^{p_1,...,p_n}$ into $A_{\alpha_1,...,\alpha_n}^{p_1,...,p_n}$ space.
\end{theorem}
%% Вставить теорию

The (weighted) Bergman projection $P_\nu$ is the orthogonal projection from the Hilbert
space $L_\nu^2(T_\Omega)$ onto its closed subspace $A_\nu^2(T_\Omega)$ and it is given by the following integral
formula
$$ P_\nu f(z)=C_\nu\int\limits_{T_\Omega} B_\nu(z,\omega)f(\omega)dV_\nu(\omega),$$
where
$B_\nu(z,\omega)= C_\nu\Delta^{-\nu+\frac{n}{r}}(\frac{z-\bar{\omega}}{i})$ is the Bergman reproducing kernel for $A_\nu(T_\Omega).$\\
Here we used the notation $dV_\nu(\omega)=\Delta^{\nu-\frac{n}{r}}(\nu)dud\nu.$ We denote by $dV(\omega)$ or $d\nu(\omega)$ the Lebegues measure on tubular domain over symmetric cone. Below and here we use
constantly the following notations $\omega=u+i\nu \in T_\Omega$ and also $z=x+iy\in T_\Omega.$ For any $f$
function from $A_\tau^{\infty}$ for large enough $\nu$ we have
$$ f(z)=C_\nu \int\limits_{T_\Omega}B_\nu(z,\omega)f(\omega)d\nu(\omega).$$
We assume that for weighted Bergman kernel $(K_\tau(z,\omega))$ the following estimate
$$\int\limits_{B(z,\tau)} (\delta(\tilde{z})^\alpha |K_\tau(\tilde{z},\omega)|dV(\tilde{z})\leq C|K_{\tau+n+1}(\omega,z)|\delta(z)^\alpha; \omega,z\in D $$
is valid for $\alpha>0, \tau\geq 0.$

Let $X$ is Herz type spaces in pseudo convex or tube domains. Let $f\in H(D_1),$ where $D_1$ is a unit disk, let $p\leqq\leq 1$ then we denote
$$Y=\tilde{X}_{p,q,\alpha,\alpha_1}=\tilde{X}(D_1) \{f\in H(D_1): \int\limits_{D_1}( \int\limits_{D(z,r)} |f(\omega)|^p dV_{\alpha}(\omega))^{frac{q}{p}}dV_{\alpha_1(z)}<\infty \}, $$
where $H(D_1)$ is a space of all analytic functions in $D_1$, $dV(\omega)$ is a Lebegues measure
on $D_1; dV_{\alpha}(z)=(1-|z|)^{\alpha}dV(z), \alpha>-1$ and $D(z,r)$ is the Bergman ball in the unit disk.

$q\leq1$, $\alpha_1>-1$,
Note very similarly we can define this analytic Herz spaces in tubular and bounded strongly pseudoconvex domains for same values of parameters using definitions of Bergman ball in tubular and bounded strongly pseudoconvex domains. These are obviously direct generalazation of classical Bergman spaces in mentioned domains.
Obviously this Herz type space can be defined also for all positive values of $p$ and $q$.
In the following theorem obtained in the paper of Shamoyan and Shipka we provide a new result on Bergman projection in these spaces under certain additional condition on Bergman kernel which we mentioned above generalizing classical Bergman projection theorem in mentioned domains. This result had also interesting applications in mentioned paper.
This new interesting result on Bergman projection relates spaces with the same dimension hovewer we provided this result for completeness of our exposition.
We refer for definition of Bergman type integral operator in this theorem to mentioned paper.

This result with very similar proof may be valid also in many other difficult domains which we mentioned above. We denote in theorem 11 below by $X$, $\tilde{X}$ mentioned Herz spaces (omitting indexes) in bounded pseudoconvex and tubular domains.

\begin{theorem}\cite{4a} %9 Теорема 3, n.22
The $P_\beta^+$ integral operator (the Bergman projection with positive kernel) for $\beta>\beta_0,$ where $\beta_0$ is large enough is mapping from $X$ to $X$ and from $\tilde{X}$ to $\tilde{X}$ for all $q\leq p\leq 1, \alpha>-1, \alpha_1>1, \tilde{\alpha}\geq 0.$
\end{theorem}

In two theorems below first for $p>1$ then for positive other values of $p$, we prove the boundedness of $S_{a,b}$ Bergman type integral operators between BMOA type spaces in the unit ball and polyball (spaces with different dimension).
These results (with not very difficult proofs based on so called composition formula) may be valid also with the same proof in Siegel domains of second type, bounded symmetric domains, and various matrix domains and in the polydisk also. Appropriate modification of $S_{a,b}$ first must be defined.
We pose this easy but intresting task to interested readers.
These results of theorem 12, 13 were used by first author many year ago to get sharp traces of analytic BMOA type spaces in the polyball.

\begin{theorem}\cite{5} %10, n.23
Let $1<p<\infty.$ Suppose $s_1,...,s_m>-1$ and $r_1,...,r_m\geq 0$ are
such that for each $j=1,...,m,$ we have $-pa_j<\min \{s_j+1, s_j+1+n-r_j\}$ and $ms_j+1< p(mb_j-n) -(m-1)(n+1).$ Denote $t=(m-1)(n+1)+\sum_{j=1}^m s_j.$ Then there is a constant $C > 0$ such that
$$\int\limits_{B_n}...\int\limits_{B_n}|S_{a,b}f(z_1,...,z_m) |^p \prod\limits_{j=1}^m \frac{(1-|z_j|^2)^{s_j}}{|1-\langle u_j,z_j\rangle|^{r_j}}d\nu(z_1)...d\nu(z_m)\leq C\int\limits_{B_n}|f(w)|^p \frac{|1-|w|^2|^t}{\prod_{j=1}^m (1-\langle u_j,w\rangle)^{r_j}}d\nu(w),
$$
for all $f\in L^p(B_n, d\nu_t)$ and $u_1,...,u_m\in B_n.$
\end{theorem}
For the case $0 < p \leq 1$ we have the following result.
\begin{theorem}\cite{5} %11, n.24
Let $0<p\leq 1.$ Suppose $s_1,...s_m>-1$ and $r_1,...,r_m\geq 0$ are
such that for each $j=1,...,m$ we have $-pa_j< \min \{s_j+1, s_j+1+n-r_j\}$ and $s_j+1<pb_j-n.$ Denote $t=(m-1)(n+1)+\sum\limits_{j=1}^m s_j. $ Then there is a constant $C>0$ such that
$$\int \limits_{B_n}...\int\limits_{B_n} |(S_{a,b}f)(z_1,...,z_m)|^p \prod\limits_{j=1}^m \frac{(1-|z_j|^2)^{s_j}}{|1-\langle u_j,z_j\rangle|^{r_j}}d\nu(z_1)...d\nu(z_m)\leq C \int\limits_{B_n}|f(w)|^p\frac{(1-|w|^2)^t}{\prod\limits_{j=1}^m|1-\langle u_j,w\rangle|^{r_j}}d\nu(w), $$
for all $f\in A^p(B_n,d\nu_t)$ and $u_1,...,u_m\in B_n.$
\end{theorem}
We need the following estimates for our results (see\cite{d3}).
\begin{equation}\label{9}
\int\limits_{\Omega}\delta^t(z)K_{n+1+t+s}(z,w)K_r(z,w)d\nu(z) \leq c_1\delta^{-s}K_r(w,\nu),
\end{equation}
where $w,\nu\in \Omega, K_t$ kernel is a function defined via estimate of Henkin -Ramirez $\phi$ function
$$|K_t(z,w)|\leq c|\phi(z,w)|^{-1}, z,w\in\Omega,$$
see \cite{d4}. In case of tubular domain for $B_{\alpha}(z,w)$ Bergman kernel we need the following estimate
\begin{equation}\label{10}
\int_{T_\Omega}\frac{\Delta^t(Im (w))d\tilde{V}(w)}{\Delta \big(\frac{z-\bar{w}}{i}\big)^{\frac{2n}{\tau}+t+s}\Delta\big(\frac{\nu-w}{i}\big)^r}\leq c_2\frac{\Delta^{-s}(Im(w))}{\Delta^r(\frac{w-\nu}{i})}, w,\nu\in T_\Omega.
\end{equation}
Here $c_2$ and $c_1$ are constants, and $t>-1, s>0, 0\leq r\leq n+1+t$ in pseudoconvex
domains and $t>-1, s>0, 0\leq r< \frac{2n}{\tau}+t$ in tubular domain, where $\tau$ is a rank of
our tube domain.

We always assume $K_t,t\in \mathbb{N}$ below, so
we consider Bergman Kernel in pseudoconvex domains only with natural index. We
define the BMOA type space in products of bounded pseudoconvex domains as a subspace of $H(\Omega^m)$ with the following finite quasinorm,for positive values of parameters involved, for all $s_j>-1, r_j\in \mathbb{N}$ for all $j$
$$\int\limits_{\Omega}...\int\limits_{\Omega} |(f)(z_1,...,z_m)|^p \prod\limits_{j=1}^m \delta^{s_j}(z_j) \prod\limits_{j=1}^m |K_{r_j}(z_j,u_j)|d\nu(z_1)...d\nu(z_m).$$
Putting $m = 1$ in this quasinorm we get such type new analytic function spaces in
pseudoconvex domains $\Omega$. We define the following Bergman type projections in
pseudoconvex domains,for large enough parameters involved.
$$S_{a,b}(f)(z_1,...,z_m)=\prod\limits_{j=1}^m\delta^{a_j}(z_j) \int\limits_{\Omega}f(w)\prod\limits_{j=1}^m |K_{a_j+b_j}(z_j,w)|(\delta(w))^{-n-1+\sum\limits_{j=1}^m b_j} d\nu(w).$$

We in our theorems 14, 15 provide under certain additional  conditions (1) and (2) on Bergman kernel new results on boudedness of $S_{a,b}$ Bergman type integral operators in so called new BMOA type  function spaces of analytic functions in tubular and bounded strongly pseudoconvex domains. 
We remark that the first author earlier proved this result also in context of the unit ball without any additional condition on Bergman kernel.
In all domains proofs of these results are similar.
Note that complete analogues of composotion formulas (1) and (2) in $\mathbb{R}^n$ are also valid. These facts are well-known to experts.

\begin{theorem}\cite{5} %12, n.25
Let $1<p<\infty, s_j>-1, r_j\in \mathbb{N}, a_j>a_0, b_j>b_0, a_0=a_0(s_1,...,s_m,p,m,n), b_0=b_0(s_1,...,s_m,p,m,n), j=1,...,m.$\\
Let $t=(m-1)(n+1)+\sum\limits_{j=1}^{m}s_j $ and (\ref{9}) holds then there is a constant $C>0$ such that
$$ \int\limits_{\Omega}...\int\limits_{\Omega}|S_{a,b}(f)(z_1,...,z_m)|^p \prod\limits_{j=1}^m \delta^{s_j}(z_j) \prod\limits_{j=1}^m|K_{r_j}(z_j,u_j)|d\nu(z_1)...d\nu(z_m)\leq C \int\limits_{\Omega} (f(w))^p \delta^t(w)\prod\limits_{j=1}^m |K_{r_j}(u_j,w)d\nu(w)|$$
the conclusion of previous part of our theorem is valid also for $0<p\leq 1,$ but for
other values of $a_0,b_0$ parameters we assume in addition to (\ref{9}) that properity (C) is
valid and $\frac{r_j}{p}\in \mathbb{N}.$
\end{theorem}
We define the analytic BMOA type space in products of tubular domains over
symmetric cones as a subspace of $H(T_{\Omega}^m)$ with the following finite quasinorm, for
positive values of parameters involved, and for all $s_j>-1$ for all $j.$
$$\int\limits_{T_\Omega}...\int\limits_{T_\Omega} |(f)(z_1,...,z_m)|^p \prod\limits_{j=1}^m \Delta^{s_j}(Im(z_j))\prod\limits_{j=1}^m \Delta^{r_j} \big(\frac{\bar{z}_j-u_j}{i}\big)dV(z_1)...dV(z_m).$$
Putting $m = 1$ in this quasinorm we get such type spaces in tubular domains over
symmetric cones $T_\Omega$. We define new Bergman type integral operators in tubular
domains as follows:
$$ S_{a,b}(f)(z_1,...,z_m)=\prod\limits_{j=1}^m \Delta^{a_j} (Im(z_j))\int\limits_{T_\Omega} f(w)\prod\limits_{j=1}^m \Delta^{-(a_j+b_j)}\big(\frac{z_j-w_j}{i}\big)[\Delta(Im(w))]^{-\frac{2n}{r}+\sum\limits_{j=1}^m b_j}dV(w),$$
$z_j\in T_\Omega, j=1,...,m$ we have that, for large enough parameters involved.

\begin{theorem}\cite{5} %13, n.26
Let $1<p<\infty, s_j>-1, r_j\geq 0, a_j>a_0, b_j>b_0, a_0=a_0(s_1,...,s_m,p,m,n), b_0=b_0(s_1,...,s_m,p,m,n), j=1,...,m.$ Let $t=(m-1)\frac{2n}{r}+\sum\limits_{j=1}^m s_j$ and (\ref{10}) holds then
there is a constant $C > 0$ so that
$$\int\limits_{T_\Omega}...\int\limits_{T_\Omega} |S_{a,b}(f)(z_1,...,z_m)|^p \prod\limits_{j=1}^m \Delta^{s_j}(Im(z_j)) \prod\limits_{j=1}^m \Delta^{-r_j}\big(\frac{\bar{z}_j-u_j}{i}\big)dV(z_1)...dV(z_m) \leq $$
$$\leq C \int\limits_{T_\Omega}|f(w)|^p \Delta^t (Im(w)) \prod\limits_{j=1}^m \Delta^{-r_j}\big(\frac{u_j-\bar{w}}{i}\big)dV(w)$$
where $u_j\in T_\Omega, j=1,...,m.$
\end{theorem}

These general results may be valid also for all $p\leq1$ and also with the same proof in many other complicated domains.

In the following theorem we extend already known result obtained by B. Sehba in classical Bergman spaces in tube in two directions simultaneously, namely to product domains and mixed norm function spaces. Earlier this result was obtained in very simple case of unit disk by O. Yaroslavceva. Proofs in both cases are similar.
In $\mathbb{R}^n$ these classes were studied by Benedek and Panzone much earlier.  Further in theorems 17, 18 we extend the classical notion of Bergman integral operator in another way defining $\mathbb{R}_X$ integral operators on product domains. They were studied in unit ball by the first author earlier. Here we managed to generalize known results from the unit disk polydisk and unit ball to more general domains, all these three new theorems are also valid in context of tube and bounded strongly pseudoconvex domains simultanuausly with the same proof.

\begin{theorem} \cite{1} %1 n.14
Let
$$ T_{\overrightarrow{\beta}}f(\overrightarrow{z})=\int\limits_{T_\Omega^m}\frac{f(w_1,...,w_m)\prod\limits_{j=1}^m \Delta^{\beta_j-\frac{n}{r}}(w_j)dv(w_j)}{\Delta^{\beta_1+\frac{n}{r}}\big ( \frac{z_1-\overline{w}_1}{i} \big)...\Delta^{\beta_m+\frac{n}{r}}\big( \frac{z_m-\overline{w}_m}{i}\big)},$$
$dv(w)=dudv, w=u+iv\in T_\Omega, \overrightarrow{z}=(z_1,...,z_m)\in T_\Omega.$

Let $\beta_j>\beta_0, j=1,...,m,$ for some fixed enough large $\beta_0.$ Then $T_{\overrightarrow{\beta}}$ operator maps $L^{\overrightarrow{p}}_{\overrightarrow{\nu}}(T_\Omega^m) $ into $A^{\overrightarrow{p}}_{\overrightarrow{\nu}}(T_\Omega^m), p_j>1, \nu_j>\frac{n}{r}-1, j=1,...,m. $
\end{theorem}

\begin{theorem} \cite{1} %2, n.15
Let $s_j>-1$ and $ms_j+1> m(\frac{2n}{r}-y_j)-(m-1)(\frac{2n}{r}), j=1,...,m.$
 Then there is exist a constant $C>0$ such that
 $$ \int\limits_{T_\Omega}|R_{x,y}g)(w)|\cdot \Delta(Im w)^{(m-1)\frac{2n}{r}+\sum\limits_{j=1}^{m}s_j}dV(w)\leq
 C\int\limits_{T_\Omega}...\int\limits_{T_\Omega}g(z_1,...,z_m)\cdot \prod\limits_{j=1}^m (\Delta^{s_j}(Im z_j))dV(z_j).$$
\end{theorem}

\begin{theorem}\cite{2} %3 (теорема 1), n.16
Let
$$[R_{x,y}(g)(w)]=
\big [ (Im w)^{-m\frac{2n}{r}+\sum\limits_{j=1}^m y_j} \big ]\int\limits_{T_\Omega}...\int\limits_{T_\Omega} g(z_1,...,z_m)\frac{[\prod\limits_{j=1}^m (Im z_j)^{x_j}]dv(z_1)...dv(z_m)}{\prod\limits_{j=1}^m \big | \Delta^{x_j+y_j}(\frac{w-z_j}{i}) \big |}, $$
for $g\in L^1(T_\Omega^m, dv(z_1),...,dv(z_m)), w\in T_\Omega.$ Let $1< p<\infty, s_j> (\frac{n}{r}-1), j=1,...,m. $ Let also for each $j=1,...,m, x_j>x_0, y_j>y_0,$ where $x_0,y_0$
are large enough positive numbers. Than there is a constant $C>0$ such
that
$$\int\limits_{T_\Omega}|(R_{x,y}(w))|^p (Im w)^{(m-1)(\frac{2n}{r})+\sum\limits_{j=1}^m (s_j-\frac{n}{r})}dv(w)\leq c\int\limits_{T_\Omega}...\int\limits_{T_\Omega} |g(z_1,...z_m)|^p \prod \limits_{j=1}^m (Im z_j)^{s_j-\frac{n}{r}}dv(z_j).$$
\end{theorem}

Let $\vec{\alpha}=(\alpha_1,...,\alpha_m), \vec{\beta}=(\beta_1,...,\beta_m)$
or $\vec{\beta}=(\beta,...,\beta).$

We modify the $R_{x,y}$ operator we defined above. Let also for $g\in L^1(T_\Omega^m)$
$$ [G_{\vec{\alpha},\beta}(g)](x+iy_1,...,x+iy_m)=\int\limits_{T_\Omega}\frac{g(w)[\Delta^{\beta}(Im w)]dv(w)}{\bigg[\prod\limits_{j=1}^m \Delta^{\alpha_j}(\frac{w-(\bar{x}+\bar{y}_j)}{i})\bigg]};$$
$x\in \mathbb{R}^n, y_j\in \Omega, j=1,...,m, \alpha_j>0, \beta>\frac{n}{r}-1, j=1,...,m.$
$$ [G_{\vec{\alpha},\vec{\beta}}(g)](x_1+iy,...,x_m+iy)=\int\limits_{T_\Omega}\frac{g(w)[\Delta^{\beta}(Im w)]dv(w)}{\bigg[\prod\limits_{j=1}^m \Delta^{\alpha_j}(\frac{w-(\bar{x_j}+i\bar{y}_j)}{i})\bigg]};$$
$\beta>\frac{n}{r}-1, j=1,...m, \alpha_j>0, x_j\in \mathbb{R}^n, j=1,...,m, y\in \Omega.$

\overfullrule=0pt In the following theorem 19 various new estimates for various new Bergman type integral operatots are provided in semiproducts and products of tube domains. They extend known results from the unit disk. They may have  also complete analogues in bounded strongly pseudoconvex domains. We pose this as a problem.
We refer for definitions of these operators to papers of first author with S. Kurilenko and related papers which are indicated in references.
One of these Bergman type integral operators ($G_{x,y}$ operators) was defined above.

\begin{theorem} \cite{2} % 4 (теорема 3), n.17
For $1\leq p<\infty$ some $ \alpha_j\in (\alpha_0,\alpha_0^{'}); \beta_j\in(\beta_0, \beta_0^{'}); \beta\in(\tilde{\beta}_0, \tilde{\beta}_0^{'}), \nu_j>(\frac{n}{r}-1), j=1,...,m,$ for some fixed positive $\alpha_0^j,(\alpha_0^j)^{'},\beta_0^j, (\beta_0^j)^{'},\tilde{\beta}_0^j, (\tilde{\beta}_0^j)^{'}, j=1,...,m. $

The following estimates are valid\\
1). $\|G_{\vec{\alpha},\beta}(g)\|_{(A_{\vec{\nu}}^p)_3} \leq c_1\|g\|_{(A_\tau^p)(T_\Omega)};$ for some values $\vec{\nu}$ and $\tau.$\\
2). $\|G_{\vec{\alpha},\beta}(g)\|_{(A_{\nu}^p)_2} \leq c_2\|g\|_{(A_\tau^p)(T_\Omega)};$ for some values $\vec{\nu}$ and $\tau.$\\
3). $ \|V_{\vec{\alpha},\vec{\beta}}(g)\|_{(A_\nu^p)_1}\leq c_3\| g\|_{(A_{\vec{\tau}}^p)_3(T_\Omega)};$ for some values $\nu$ and $\vec{\tau}.$\\
4). $ \|U_{\vec{\alpha},\vec{\beta}}(g)\|_{(A_\nu^p)_1}\leq c_4\| g\|_{(A_{\vec{\tau}}^p)_2(T_\Omega)};$ for some values $\nu$ and $\vec{\tau},$\\
where $\alpha_0^j,...,(\beta_0^j)^{'}$ depend on $\nu_j, \tau_j,p,n,\nu, \tau, j=1,...,m$ and $\frac{1}{p}+\frac{1}{p'}=1.$
\end{theorem}

%Тут был кусок текста, теперь он в конце первого параграфа

In the following theorem we formulate a result which relates spaces with different dimensions in tube and pseudoconvex domains. 
We would like to remark for $m=1$ case the additional condition in theorem 20 vanishes and we get a known result in mentioned domains
This theorem we formulate in context of tube and bounded pseudoconvex domains noting that arguments of proofs were taken from unit disk and polydisk case which was provided by other authors earlier. Nevertheless in these general domains this result is interesting enough.

 \begin{theorem} \cite{3} %5 n.18
 
 1). Let $1<p_j<\infty, j=1,...,m,$ and let
 $$T_\beta f(\overrightarrow{z})=\int\limits_{T_\Omega}\frac{f(w)d\nu_{\beta_1}(w)}{\prod\limits_{j=1}^m \Delta^{\beta+\frac{2n}{r}}\big(\frac{z_j-w}{i}\big)},$$
 where $\overrightarrow{z}=(z_1,...,z_m), z_j\in T_\Omega, j=1,...,m.$

 Let $\beta>\beta_0, j=1,...,m,$ for some fixed large enough $\beta_0.$ Then $T_\beta$ maps $A_\nu^{p_m}(T_\Omega)$ or $L_\nu^{p_m}$ to $A_{\overrightarrow{\nu}^{\overrightarrow{p}}}(T_\Omega^m),$
 $$\nu=\sum\limits_{j=1}^{m-1} \big[\nu_j+\frac{2n}{r}\big]\big(\frac{p_m}{p_j}\big)+\nu_m,$$
 here $A_\nu^p=A_\nu^P(T_\Omega^1),$ if for $F\in A^{\overrightarrow{q}}(\overrightarrow{\alpha})(T_\Omega^m) F(w,...,w)\in A^{q_m}(\alpha)(T_\Omega)$ where
 $$ \alpha=\beta-\frac{\tau q_m}{p_m}, \tau=\alpha_m+\sum\limits_{j=1}^{m-1}\big( \alpha_j+\frac{2n}{r}\big)\frac{p_m}{p_j},$$
 $$\beta_1=(\beta+\frac{2n}{r})m-\frac{2n}{r}, \frac{1}{p_j}+\frac{q}{q_j}=1, j=1,...,m.$$
  2).Let
  $$(S_{\beta_0})(f(\overrightarrow{z}))=\int\limits_D f(w)\prod\limits_{j=1}^m K_{\tilde{\beta}_0}(z_j,w)d\nu_{\beta_1}(w),$$
  $\overrightarrow{z}=(z_1,...,z_m), z_j\in \Lambda, j=1,...,m.$

  Let $\beta>\tilde{\beta}^*_0, j=1,...,m,$ for some fixed large enough $\tilde{\beta}^*_0, 0<p_j<\infty.$
  Then $(S_{\overrightarrow{\beta}})$ operator maps $A_\nu^{p_m}(\Lambda)$ or $L_\nu^{p_m}$ into $A_{\overrightarrow{\nu}}^{\overrightarrow{p}}(\Lambda^m),$
  $$\nu=\sum\limits_{j=1}^{m-1}[\nu_j+(n-1)]\big(\frac{p_m}{p_j}\big)+\nu_m, 1<p_j<\infty, \nu_j>-1, j=1,...,m, \nu>-1,$$
  if for $F\in A_{\overrightarrow{\alpha}}^{\overrightarrow{q}}(\Lambda^m) F(w,...,w)\in A^{q_m}(\tilde{\alpha})(\Lambda),$ where
  $$\tilde{\alpha}=\beta_0-\frac{\tau q_m}{p_m}, \tau =\alpha_m+\sum\limits_{j=1}^{m-1} (\alpha_j+n+1)\frac{p_m}{p_j},$$
  $\tilde{\beta}_0=\beta_0+n+1, \beta_1=(\beta_0+n+1)m-(n+1), \frac{1}{p_j}+\frac{1}{q_j}=1, j=1,...,m.$

 \end{theorem}

%Здесь был кусок текста, теперь он в конце первого параграфа

In the following two theorems we show the boundedness of certain Herz type integral operators (similar to Bergman type operators) in some function spaces with different dimensions in tube domains in $\mathbb{C}^n$.
Proofs of these theorems are purely technical and they were provided in papers of the first author with S. Kurilenko.
We refer the reader to function spaces we mention in those theorems to papers we mention in those theorems. These results are new as far as we know even in case of simplest unit disk and polydisk.

\begin{theorem} \cite{a} %n.1
For $1\leq p<\infty$ the following estimate is valid
$$ \|T_{\vec{\alpha},\beta,\gamma}(g)\|_{(L_{\vec{\nu}}^p)_1(T_\Omega^m)}\leq c\|g\|_{(L_\tau^p)(T_\Omega)},$$ where $\frac{1}{p}+\frac{1}{p'}=1, \beta>-1, p'\gamma>-1,$ and
$$ \alpha_j>\frac{\beta}{m}+\frac{2n}{rm}+\frac{3n}{rp}-\frac{1}{p}+\frac{3n}{rp'm}-\frac{1}{p'm}+\frac{\gamma}{m}, j=1,...,m,$$ and
$$ \tau=p\beta+\frac{2nm}{r}-p\sum\limits_{j=1}^m (\alpha_j-\frac{\nu_j}{p})+ p\gamma+ \frac{4np}{rp'}-\frac{nm}{rp}.$$
\end{theorem}
\begin{theorem} \cite{a} %n.2
For $1\leq p<\infty$ the following estimate is valid
$$\|T_{\vec{\alpha,\beta,\gamma}}(g)\|_{(L_{\vec{\nu},\vec{s}}^p)_2 (T_\Omega^m)}\leq c\| |g| \|^p_{(L_\tau^p)(T_\Omega)},$$
where $\frac{1}{p}+\frac{1}{p'}=1, \beta>-1, p'\gamma>-1, s_j>\frac{n}{r}-1,$
$$ \alpha_j>\frac{\beta}{m}+\frac{2n}{rm}+\frac{2n}{rp}-\frac{1}{p}+\frac{3n}{rp'm}-\frac{1}{p'm}+\frac{\gamma}{m}+\frac{s_j}{p}, j=1,...,m,$$
and
$$ \tau=p\beta+\frac{2nm}{r}+\sum\limits_{j=1}^m (s_j-p\alpha_j +\nu_j) +p\gamma+\frac{4np}{rp'}.$$
\end{theorem}
Let $\alpha_j>0, b_j>-1, \gamma_j>-1, j=1,...,m.$ m. Now we define another new
integral Herz type operator
$$ [T^1_{\vec{\alpha},\vec{\beta},\vec{\gamma}}(g)](z_1,...,z_m)=$$
$$=\int\limits_{T_\Omega}...\int\limits_{T_\Omega} \int\limits_{B(\tilde{w}_1,r)}...\int\limits_{B(\tilde{w}_m,r)} g(w_1,...,w_m)\frac{\prod_{j=1}^{m} [\Delta^{\beta_j}(Im w_j)]dv(w_1)...dv(w_m)}{\big [\prod_{j=1}^m \Delta^{\alpha_j}(\frac{z_j-\bar{w}_j}{i})\big] } \prod\limits_{j=1}^m \Delta^{\gamma_j} (Im \tilde{w}_j)dv(\tilde{w}_1)...dv(\tilde{w}_m),$$
where $ z_j\in T_\Omega, j=1,...,m$ for a $g$ function from $L^1$
class on product of tubes. These type Bergman integral operators are new even in one dimensional case.

In the following four theorems we provide new results on new Bergman type integral operators 
In tube domains obtained in papers of the first author with S. Kurilenko.
These type results may be also valid in bounded strongly pseudoconvex domains.
These new assertions contain long technical calculations but they are new even in one dimensional case. These new Bergman type (Herz type) integral operators are even new in the case of unit disk.

\begin{theorem} \cite{a} %n.3
For $1\leq p<\infty$ the following estimate is valid
$$\| T^1_{\vec{\alpha},\vec{\beta},\vec{\gamma}}(g)\|_{(L_{\vec{\nu}}^p)_1(T_\Omega^m)}\leq c\| g\|_{(L_\tau^p)(T_\Omega^m)},$$
where $\frac{1}{p}+\frac{1}{p'}=1, \beta_j>-1, p'\gamma_j>-1,$
$$\alpha_j>\beta_j+\frac{2n}{rp'}+\frac{3n}{r}+\gamma_j-1, j=1,...,m, $$
and
$$ \tau_j=p\beta_j +\frac{2n}{r}-p\alpha_j+\nu_j+\gamma_j p+\frac{4np}{rp'}+\frac{n}{r}, j=1,...,m.$$
\end{theorem}
\begin{theorem} \cite{a} % n.4
For $1\leq p<\infty$ the following estimate is valid
$$\| T^1_{\vec{\alpha},\vec{\beta},\vec{\gamma}}(g)\|_{(L_{\vec{\nu}}^p)_2(T_\Omega^m)}\leq c\| g\|_{(L_\tau^p)(T_\Omega^m)},$$
where $\frac{1}{p}+\frac{1}{p'}=1, \beta_j>-1, p'\gamma_j>-1,$
$$\alpha_j>\beta_j+\gamma_j+\frac{s_j}{p}-1+\frac{3n}{r}+\frac{n}{p}+\frac{2n}{rp'}, j=1,...,m,$$
and
$$ \tau_j=p\beta_j +\frac{2np}{r}+s_j-p\alpha_j +\nu_j+p\gamma_j+\frac{2np}{rp'}, j=1,...,m.$$
\end{theorem}

%The weighted Bergman kernel $B_\nu$ of $T_\Omega$ is given as usual by
%$$B_\nu(w,z)=(d_\nu)\Delta \big(\frac{w-\bar{z}}{i}\big)^{-\nu-\frac{n}{r}}, w,z\in T_\Omega, \nu\in \mathbb{R},$$
%is a Bergman constant, where
%$$ d_\nu=(c_\nu^{-1})\Gamma \big (\nu+\frac{n}{r}\big).$$

%Let $\Omega$ be an irreducible symmetric cone in the Euclidean space $V$ , and $T_\Omega=V+i\Omega$
%the corresponding tube domain in the complexified space $V^{\mathbb{C}}.$\\
%Let $T_\Omega$ be the tube domain over symmetric cone, and $H(T_\Omega)$ be the space of all
%analytic functions in tube (see \cite{d1}). We define Bergman spaces for %$1\leq p,q<\infty, \gamma>\frac{n}{r}-1.$ Let
%$$(A_{\gamma}^{p,q})(T_\Omega\bigg \{ f\in H(T_\Omega): \bigg(\int\limits_\Omega \big ( \int\limits_{\mathbb{R}^n}|f(x+iy)|^p dx \big)^{\frac{q}{p}}(\Delta^{\gamma-\frac{n}{r}}(y))dy\bigg)^{\frac{1}{q}}<\infty \bigg\}.$$

%Replacing $A$ by $L$ and $H$ by $L^1$ we get as usual known larger spaces of measurable
%functions in tube $T_\Omega.$ Note $(A^{p,q}_\gamma)=\{0\}$ if $\gamma \leq \frac{n}{r}-1.$

\begin{theorem} \cite{b} % n.5
There are $\nu_1=\nu_1(\alpha,n,r,q), \nu_2=\nu_2(\alpha,n,r,q),$ so that for $1\leq p,q<\infty, \nu\in\mathbb{R}, \gamma=\alpha+\beta+\frac{n}{r}, \alpha+\beta>-1$ then $T^{+}_{\alpha,\beta,\gamma}$ is a bounded operator on $L_\nu^{p,q}(T_\Omega)$ for all $\nu\in (\nu_1,\nu_2).$
\end{theorem}
\begin{theorem} \cite{b} %n. 6
Let $(Q^+)$ be $(T^+_{\alpha,\beta,\gamma})$ operator for $\alpha=0, \gamma=\nu+m; \beta=\nu-\frac{n}{r}.$ Then $(Q^+)$ for $\nu+m>\frac{n}{r}-1, 1\leq p,q<\infty,$ is a bounded operator from $L_\nu^{p,q}$ to $L_{\nu+mq}^{p,q}$ if $\nu\in (\nu_1,\nu_2)$ for some $\nu_1=\nu_1(\alpha,n,r,q), \nu_2=\nu_2(\alpha,n,r,q), (T^+_{\alpha,\beta,\gamma})$ is a bounded
operator on $L^{\infty}$ if $\alpha>\frac{n}{r}-1, \beta>-1, \gamma=\alpha+\beta+\frac{n}{r}.$ The same is valid for $T_{\alpha,\beta,\gamma}$
operator.
\end{theorem}
Let
$$ (T_\beta h)(\vec{w})=\underbrace{\int\limits_{\mathbb{R}^n}...\int\limits_{\mathbb{R}^n}}_m \int\limits_\Omega \frac {h(x_1+iy,...,x_m+iy)\Delta^\beta (y)}{\prod\limits_{j=1}^m \Delta^{\beta_j+\frac{n/r+mn/r}{m}}(\frac{\bar{x}_j+i\bar{y}-w_j}{i})}dx dy_1...dy_m,$$
where $\vec{w}=\{w_1,...,w_m\}=\{\zeta_1+i\eta,...,\zeta_m+i\eta\}\in T_\Omega,$
$$(\tilde{T}_{\vec{\beta}}h)(\vec{w})=\int\limits_{\Omega^m}\int\limits_{\mathbb{R}^n} \frac{h(x+iy_1,...,x+iy_m) \Delta^{\beta_1}(y_1)...\Delta^{\beta_m}(y_m)}{\prod\limits_{j=1}^m \Delta^{\beta_j+\frac{n/r+mn/r}{m}}(\frac{\bar{x}+i\bar{y}_j-w_j}{i})}dx dy_1...dy_m,$$
where $\vec{w}=\{w_1,...,w_m\}=\{\zeta+i\eta_1,...,\zeta+i\eta_m\} \in T_\Omega,$ and $\vec{\beta}=\{\beta_1,...,\beta_m\}, \beta>\frac{n}{r}-1, \beta_j>\frac{n}{r}-1, j=1,...,m, h\in L^1(T^m_\Omega).$
\begin{theorem}\cite{b} % n.7
Let $1<p<\infty, \beta>m\frac{m}{r}-m-\frac{n}{r}+\frac{mnp}{r(p-1)}, \tau>(p-1)(\frac{2mn}{r}-m-\frac{n}{r}-1)-1$ and $m\in\mathbb{N}, m>1, \tau-p\beta< 1-\frac{n}{r}$ then
$$ \underbrace{\int\limits_{\mathbb{R}^n}...\int\limits_{\mathbb{R}^n}}_m \int\limits_\Omega |(T_\beta h)(\vec{w})|^p\Delta^\tau(\eta) d\eta d\zeta_1...d\zeta_m \leq c\underbrace{\int\limits_{\mathbb{R}^n}...\int\limits_{\mathbb{R}^n}}_m \int\limits_\Omega |h(\vec{z})|^p \Delta^\tau(y)dy dx_1...dx_m, $$
where $\vec{w}=\{w_1,...,w_m\}=\{\zeta_1+i\eta,...,\zeta_m+i\eta\}\in T^m_\Omega, \vec{z}=\{z_1,...,z_m\}=\{x_1+iy,...,x_m+iy\}\in T^m_\Omega.$
\end{theorem}

\begin{theorem}\cite{b} % n.8
Let $\beta_j>\frac{n}{r}-1, j=1,...,m, 1<p<\infty,$ $\beta_j>2+\frac{n}{r}-\frac{2n}{pr}-\frac{3}{m}-\frac{2}{p}+\frac{n}{rmp}+\frac{2}{mp}, \tau_j>\frac{n}{r}p-p-\frac{n}{r}, \tau_j-p\beta_j<1-\frac{n}{r}, j=1,...,m,$ then
$$\int\limits_{\Omega^m} \int\limits_{\mathbb{R}^n} |(\tilde{T}_{\vec{\beta}}h)(\vec{w})|^p \Delta^{\tau_1}(\eta_1)\times...\times\Delta^{\tau_m}(\eta_m)d\eta d\zeta_1...d\zeta_m \leq c\int\limits_{\Omega^m} \int\limits_{\mathbb{R}^n}|h(\vec{z})|^p \Delta^{\tau_1}(y_1)\times...\times \Delta^{\tau_m}(y_m) dx dy_1...dy_m,$$
where $\vec{w}=\{\zeta+i\eta_1,...,\zeta+i\eta_m\}, \vec{z}=\{z_1,...,z_m\}=\{x+iy_1,...x+iy_m\}, \vec{z}\in T^m_\Omega, \vec{w}\in T^m_\Omega.$
\end{theorem}
Next let $1\leq p<\infty, f=f(z_1,...z_m),$ we consider analytic subspaces of $H(T_\Omega^m), T_\Omega^m=T_\Omega\times...\times T_\Omega, \nu_j>\frac{n}{r}-1, \nu>\frac{n}{r}-1, j=1,...,m.$ These are spaces $(A_\nu^p)_1,(A_\nu^p)_2,(A_\nu^p)_3 $ with norms
$$ \|f\|^p_{(A_{\vec{\nu}}^p)_1}=\int\limits_{T_\Omega}....\int\limits_{T_\Omega} |f(x_1+iy_1,...,x_m+iy_m)|^p \prod\limits_{j=1}^m \Delta^{\nu_j-\frac{n}{r}}(y_j)dx_jdy_j<\infty, $$
$$ \|f\|^p_{(A_{\vec{\nu}}^p)_2}=\int\limits_{\mathbb{R}^n}....\int\limits_{\mathbb{R}^n}\int\limits_\Omega |f(x_1+iy_1,...,x_m+iy_m)|^p \Delta^{\nu-\frac{n}{r}}(y) \big( \prod\limits_{j=1}^m dx_j \big)dy<\infty , $$
$$ \|f\|^p_{(A_{\vec{\nu}}^p)_3}=\int\limits_{\mathbb{R}^n}\int\limits_\Omega....\int\limits_{\Omega} |f(x_1+iy_1,...,x_m+iy_m)|^p \prod\limits_{j=1}^m \Delta^{\nu_j-\frac{n}{r}}(y_j)dxdy_j<\infty. $$
Let also for $g\in L^1(T_\Omega^m)$
$$ (V_{\vec{\alpha},\vec{\beta}}g)(\vec{w})= \underbrace{\int\limits_\Omega...\int\limits_\Omega}_m \int\limits_{\mathbb{R}^n} \frac{g(x+iy_1,...,x+iy_m) (\Delta y_1)^{\beta_1}\times...\times(\Delta y_m)^{\beta_m}}{\prod\limits_{j=1}^m \Delta^{\alpha_j}\big(\frac{\bar{x}+i\bar{y}_j-w_j}{i}\big)}dx dy_1...dy_m,$$
$$ (U_{\vec{\alpha},\vec{\beta}}g)(\vec{w})= \underbrace{\int\limits_{\mathbb{R}^n}...\int\limits_{\mathbb{R}^n}}_m \int\limits_{\Omega} \frac{g(x_1+iy,...,x_m+iy) (\Delta y)^\beta}{\prod\limits_{j=1}^m \Delta^{\alpha_j}\big(\frac{\bar{x}+i\bar{y}_j-w_j}{i}\big)}dy dx_1...dx_m,$$
$w=(w_1,...,w_m), w_j\in T_\Omega, j=1,...,m, \beta_j>\frac{n}{r}-1, \beta>\frac{n}{r}-1, \alpha_j>0, j=1,...,m, \vec{\beta}=(\beta_1,...,\beta_m)$ or $\vec{\beta}=(\beta,...,\beta).$

Let also for $g\in L^1(T_\Omega^m)$
$$[G_{\vec{\alpha},\beta}(g)](x+iy_1,...,x+iy_m)=\int\limits_{T_\Omega} \frac{g(w)[\Delta^\beta(Im w)]dv(w)}{| \prod\limits_{j=1}^m \Delta^{\alpha_j} \big ( \frac {w-(\bar{x}+\bar{y}_j)}{i} \big ) |},$$
where $x\in \mathbb{R}^n, y_j\in \Omega, j=1,...,m, \alpha_j>0, \beta>\frac{n}{r}-1, j=1,...,m;$
$$[G_{\vec{\alpha},\vec{\beta}}(g)](x_1+iy,...,x_m+iy)=\int\limits_{T_\Omega} \frac{g(w)[\Delta^\beta(Im w)]dv(w)}{| \prod\limits_{j=1}^m \Delta^{\alpha_j} \big ( \frac {w-(\bar{x}+i\bar{y}_j)}{i} \big ) |},$$
where $x_j\in \mathbb{R}^n, j=1,...,m, y\in \Omega, \beta>\frac{n}{r}-1, \alpha_j>0, j=1,...,m.$

We define below various quazinorms in products of tubular domains and on semiproducts of such type domains.
They extend classical quazinorm of analytic Bergman $A^p_\alpha(T_\Omega)$ space in tube studied by Sehba and many other authors.
We provide in theorem 29 new results on the action of some new Bergman type integral operators in analytic spaces in tube with these quazinorms extending previously known assertions obtained by B. Sehba and his coauthors.

We refer to mentioned in theorem 29 spaces and operators which are not defined the mentioned paper in this theorem 29 with S. Kurilenko.

\begin{theorem} \cite{b} % n.9
For $1\leq p<\infty$ some $\alpha_j\in (\alpha_0, \alpha_0^{'}); \beta_j\in(\beta_0, \beta_0^{'}); \beta\in(\tilde{\beta}_0, \tilde{\beta}_0^{'}),$ $\nu_j>\frac{n}{r}-1, \tau_j>\frac{n}{r}-1, j=1,...,m$ for some fixed positive $\alpha_0^j, (\alpha_0^j)', \beta_0^j, (\beta_0^j)', j=1,...,m.$ \\
The following estimates are valid:\\
$1.\|G_{\vec{\alpha},\beta}(g)\|_{(A_{\vec{\nu}}^p)_3}\leq c_1 \|g\|_{(A_\tau^p)(T_\Omega)}; $ for some values $\vec{\nu}$ and $\tau,$ \\
$2.\|\tilde{G}_{\vec{\alpha},\beta}(g)\|_{(A_{\nu}^p)_2}\leq c_2 \|g\|_{(A_\tau^p)(T_\Omega)}; $ for some values ${\nu}$ and $\tau,$ \\
$3.\|V_{\vec{\alpha},\vec{\beta}}(g)\|_{(A_{\vec{\nu}}^p)_1}\leq c_3 \|g\|_{(A_{\vec{\tau}}^p)_3(T_\Omega)}; $ for some values $\nu$ and $\vec{\tau},$ \\
$1.\|U_{\vec{\alpha},\vec{\beta}}(g)\|_{(A_{\vec{\nu}}^p)_1}\leq c_4\|g\|_{(A_{\vec{\tau}}^p)_2(T_\Omega)}; $ for some values $\nu$ and $\vec{\tau},$ \\
where $\alpha_0^j,...,(\beta_0^j)'$ depend on $\nu_j, \tau_j, p,n,\nu, \tau, j=1,...,m$ and $\frac{1}{p}+\frac{1}{p'}=1.$
\end{theorem}

In the following two theorems we provide complete descriptions of traces of certain BMOA type spaces of several variables in tube.
Same results with very similar proofs are valid in bounded pseudoconvex domains with smooth boundary.
Complete detailed proofs can be seen in context of the unit ball in papers of the first author.
All these sharp results on traces are heavily based on certain new projection theorems on $S_{a,b}$ Bergman type operators which we defined above acting between BMOA type spaces of different dimension.

\begin{theorem} \cite{21}
Let $p>1, \tau\in (0, \infty), r_j\in \mathbb{N}, s_j>-1, j=1,...,m.$ If $ t=(m+1)(2n/r)+\sum\limits_{j=1}^m s_j,$ then for $r=\sum\limits_{j=1}^m r_j$ we have $$ Trace(M^p_{r_1,...,r_m, \tau, s_1,...,s_m}(\Omega^m))= M^p_{r,p,\tau m} (\Omega) $$
for all $n, n|r>n_0, $ where $n_0=n_0(p,\tau, r_1,...,r_m, m).$
\end{theorem}

\begin{theorem}\cite{21}
Let $p>\leq 1, \tau\in (0, \infty), r_j\in \mathbb{N}, s_j>-1, j=1,...,m.$ If $r_j|p\in \mathbb{N}, j=1,...,m$ and $t=(m-1)(2n|r)+\sum\limits_{j=1}^m s_j,$ then for $r=\sum\limits_{j=1}^m r_j$ we have
$$Trace(M^p_{r_1,...,r_m, \tau, s_1,...,s_m}(\Omega^m))= M^p_{r,p,\tau m} (\Omega)$$
for all $n, n|r>n_0, $ where $n_0=n_0(p,\tau, r_1,...,r_m, m).$
\end{theorem}

Some interesting new results on Bergman type projections between Bloch type spaces of different dimensions extending some previously known results in tubular and bounded strongly pseudoconvex domains, unit ball and polydisk  were provided earlier in paper of first author and firsth author with S. Kurilenko which we mentioned for readers in references.

In \cite{k1} was provided an extension of a classical result about Bergman projection in the unit disk to pseudoconvex domains and tubular domains over symmetric cones. These proofs in particular are based on some arguments provided earlier in [11] and [12] for a less general unit ball case. New results complement similar type recent results on Bergman type projections obtained earlier in [1-5] in pseudoconvex domains and in [6-9] in tubular domains. %на сайте только аннотация, списка литературы нет
In recent two decades by various authors various interesting results were provided concerning Bergman projection in analytic function spaces of one and several variables with Muckenhoupt weights, Bekolle weights, etc. These issues will not be considered in this paper. Another group of interesting new results concerns Bergman projection theorems in analytic spaces of one and several variables with various restrictions on boundary of a domain in $\mathbb{C}$ or $\mathbb{C}^n$, we refer the interested readers to research papers of E. Stein and his various coauthors. 
All these results may be considred as extension of old classical results on Bergman projections in classical function spaces of one and several variables in the unit disk, polydisk and in the unit ball provided earlier by various authors in previous century.

Our results may be valid also with very similar proofs in various Siegel domains of second type, matrix domains, bounded symmetric domains and various other domains in $\mathbb{C}^n$ with complicated structure.

In the list of references we collected among other papers many recent papers of various authors which can be interesting for readers concerning action of Bergman type projection in various function spaces and various domains as well as some other recent research papers of the first author which we didn't mention in this paper.

Many Bergman projection theorems acting between spaces of different dimension appeared in papers of the first author as a part of solution of various trace problems in function spaces previously.
We also refer the reader for such results to papers of the first author related with trace problem posted in various databases.

Various projection theorems in harmonic function spaces in various domains (see for example [45]) can be probably extended to certain new projection theorems between harmonic function spaces with different dimensions using arguments and approaches developed in proofs of theorems of this paper in tube and pseudoconvex domains. We pose this as a problem for interested readers.

This survey concerns only function spaces of analytic functions of several variables in various complicated domains. We pose as a new interesting problem to find analogues of Bergman type projection results acting between analytic function spaces with different dimension in so called Fock-type spaces of entire functions of one and several variables based on various well-known integral representations for them. Extending various already known results in this direction provided earlier by many authors.
As far as we know this problem is new.

\end{document}